\documentclass{emsprocart}
\usepackage{amscd,amssymb}

\newtheorem{theorem}{Theorem}[section]

\newtheorem{lemma}[theorem]{Lemma}
\newtheorem{proposition}[theorem]{Proposition}
\newtheorem{conjecture}[theorem]{Conjecture}

\theoremstyle{definition}

\newcommand\alp{\alpha}		
		
\newcommand\gam{\gamma}		\newcommand\Gam{\Gamma}
\newcommand\del{\delta}		\newcommand\Del{\Delta}

\newcommand\lam{\lambda}		\newcommand\Lam{\Lambda}

\newcommand\calE{{\mathcal{E}}}
\newcommand\calF{{\mathcal{F}}}
\newcommand\calG{{\mathcal{G}}}
\newcommand\calH{{\mathcal{H}}}

\newcommand\calJ{{\mathcal{J}}}
\newcommand\calK{{\mathcal{K}}}
\newcommand\calL{{\mathcal{L}}}
\newcommand\calM{{\mathcal{M}}}

\newcommand\calO{{\mathcal{O}}}
\newcommand\calP{{\mathcal{P}}}
\newcommand\calQ{{\mathcal{Q}}}

\newcommand\calS{{\mathcal{S}}}

\newcommand\calU{{\mathcal{U}}}

\newcommand\calW{{\mathcal{W}}}
\newcommand\calX{{\mathcal{X}}}
\newcommand\calY{{\mathcal{Y}}}

		\newcommand\bfH{{\mathbf H}}

\newcommand\QQ{\mathbb{Q}}

\newcommand\PP{\mathbb{P}}
\renewcommand\AA{\mathbb{A}}

\newcommand\FF{\mathbb{F}}
\newcommand\GG{\mathbb{G}}

\newcommand\ZZ{\mathbb{Z}}

\newcommand\CC{\mathbb{C}}
\newcommand\VV{\mathbb{V}}

	\newcommand\grg{{\mathfrak{g}}}
	\newcommand\grh{{\mathfrak{h}}}

	\newcommand\grn{{\mathfrak{n}}}

	\newcommand\grr{{\mathfrak{r}}}
	\newcommand\grs{{\mathfrak{s}}}
	\newcommand\grt{{\mathfrak{t}}}

\newcommand\sdp{\times \hskip -0.3em {\raise 0.3ex
\hbox{$\scriptscriptstyle |$}}} 

\newcommand\Aut{\operatorname{Aut}}

\newcommand\Gr{\operatorname{Gr}}

\newcommand\Hom{\operatorname {Hom}}

\newcommand\Res{\operatorname{Res}}

\newcommand\Perv{\operatorname{Perv}}

\newcommand\SL{{\rm SL}}

\newcommand\supp{\operatorname{supp}}

\newcommand\Sym{\operatorname{Sym}}

\newcommand\olam{{\overline{\lambda}}}

\newcommand\hatG{{\widehat{G}}}

\newcommand\hatT{{\widehat{T}}}

\newcommand\tilG{{\widetilde{G}}}

\newcommand\x{\times}
\newcommand\ten{\otimes}

\renewcommand{\>}{\rangle}

\newcommand{\ra}{\rangle}
\newcommand{\la}{\langle}

\newcommand\qlb{{\overline \QQ}_l}

\newcommand\Bun{\operatorname{Bun}}

\newcommand\Rep{\operatorname{Rep}}

\newcommand\Vect{\VV ect}

\newcommand\Pic{\operatorname{Pic}}

\newcommand\oGr{\overline \Gr}

\newcommand\aff{\operatorname{aff}}

\newcommand\IC{\operatorname{IC}}
\newcommand\gr{\operatorname{gr}}
\newcommand\ocalW{\overline{\calW}}
\newcommand\IH{\operatorname{IH}}
\renewcommand\SL{\operatorname{SL}}
\newcommand\fsl{\mathfrak{sl}}

\newcommand\nc{\newcommand}
\nc\bln{\Bun_G(\calL^0)}
\nc\ocF{\overline \calF}
\nc\ocL{\overline \calL}
\nc\oEis{\overline \Eis}
\nc\tilBun{\widetilde \Bun}
\nc\tcF{{\widetilde \calF}}
\nc\Eis{\operatorname{Eis}}
\nc\chn{\grn^{\vee}}

\nc\chh{\grh^{\vee}}
\newcommand\kk{\Bbbk}
\nc\sph{\operatorname{sph}}

\newcommand\hatS{\widehat{S}}

\newcommand\tBun{\widetilde{\Bun}}

\newcommand\Fr{\operatorname{Fr}}
\newcommand\Sp{\operatorname{Sp}}
\renewcommand\Re{\operatorname{Re}}

\contact[braval@math.brown.edu]{Alexander Braverman, Department of Mathematics, Brown University}
\contact[kazhdan@math.huji.ac.il]{David Kazhdan, Einstein Institute of Mathematics, Hebrew University of Jerusalem}
\title[Representations of affine Kac-Moody groups over local and
global fields]{Representations of affine Kac-Moody groups over local and
global fields: a survey of some recent results}
\author[Alexander Braverman and David Kazhdan]{Alexander Braverman and David Kazhdan}

\begin{document}
\begin{abstract}
Let $G$ be a reductive algebraic group over a local field $\mathcal K$ or a global field $F$. It is well known that
there exists a non-trivial and interesting representation theory of the group $G(\mathcal K)$ as well as
the theory of automorphic forms on the adelic group $G({\mathbb A}_F)$.
The purpose of this talk is to give a survey of some recent constructions and results, which show that
there should exist an analog of the above theories  in the case when
$G$ is replaced by the corresponding affine Kac-Moody group $\hatG_{\operatorname{aff}}$ (which is
essentially built from the formal loop group $G((t))$ of $G$). Specifically we discuss the
following topics :
affine (classical and geometric) Satake isomorphism, Iwahori-Hecke algebra of $\hatG_{\aff}$,
affine Eisenstein series and Tamagawa measure.
\end{abstract}
\section{Introduction}
\subsection{Reductive groups: notations}
Let $G$ be a split connected reductive algebraic group defined over integers. We fix  a maximal split torus
$T$ in $G$ and denote by $\Lam$ its lattice of cocharacters (the coweight lattice of $G$). We let $W$ denote the Weyl
group of $G$; it acts on $T$ and $\Lam$.

Let $T^{\vee}$ denote the dual torus of $T$;   considered  as a group over $\CC$. The  lattice of characters of $T^{\vee}$ is canonically isomorphic to $\Lam$.
Let $G^{\vee}$ denote the Langlands
dual group of $G$ which by the definition contains $T^{\vee}$ as a maximal torus.
\subsection{The group $G_{\aff}$}To a split connected reductive group $G$ as above one can associate the corresponding
affine Kac-Moody group $G_{\aff}$ in the following way.

Let $\Lam$ denote the coweight lattice of $G$ let $Q$ be an integral, even,
negative-definite symmetric bilinear form on $\Lam$
which is invariant under the Weyl group of $G$.

One can consider the polynomial loop group $G[t,t^{-1}]$ (this is an infinite-dimensional group ind-scheme).
It is well-known that a form $Q$ as above gives rise to a central extension $\tilG$ of $G[t,t^{-1}]$:
$$
1\to \GG_m\to \tilG\to G[t,t^{-1}]\to 1
$$
Moreover,\ $\tilG$ has again a natural structure of a group ind-scheme.

The multiplicative group $\GG_m$ acts naturally on $G[t,t^{-1}]$ and this action lifts to $\tilG$.
We denote the corresponding semi-direct product by $G_{\aff}$; we also let $\grg_{\aff}$ denote its Lie algebra.
Thus if $G$ is semi-simple then $\grg_{\aff}$ is an untwisted affine Kac-Moody Lie algebra in the sense of \cite{Kac};
in particular,
it can be described by the corresponding affine root system.

Similarly, one can consider the corresponding completed affine Kac-Moody group $\hatG_{\aff}$ by replacing the polynomial
loop group $G[t,t^{-1}]$ with the formal loop group $G((t))$ in the above definitions.

We shall also denote by $G'_{\aff}$ (resp. $\hatG'_{\aff}$) the quotient of $G_{\aff}$ (resp. of $\hatG_{\aff}$)
by the central $\GG_m$.
\subsection{The dream}Let $\calK$ be  a non-archimedian local field
with ring of integers $\calO$ and residue field $\kk$.  A smooth representation of $G(\calK)$ is
 a vector space $V$ over $\CC$ together with a homomorphism $\pi:G(\calK)\to \Aut(V)$ such that the stabilizer
of every $v\in V$ contains an open compact subgroup $K$ of $G(\calK)$.  We denote the category of  smooth representations by
$\calM(G(\calK))$. The   category $\calM(G(\calK))$
has been extensively studied in the past 50 years. Similarly,   given a global field
$F$ we can consider automorphic representations of $G(\AA_F)$ where  $\AA_F$  is the
ring of adeles of $F$. In both (local and global) cases
the most interesting statement about the above representations is \textit{Langlands correspondence} which relates
representations of $G(\calK)$ (resp. automorphic representations of $G(\AA_F)$) to homomorphisms from the
absolute Galois group of $\calK$ (resp. of $F$) to the Langlads dual group $G^{\vee}$ of $G$.

Our dream would be to develop an analog of the above representation theories and the Langlands correspondence for
the group $G_{\aff}$ or $\hatG_{\aff}$ (or, more generally, for any symmetrizable Kac-Moody group). This
 is a  fascinating task by itself but we also believe that a fully developed theory of automorphic forms
for $\hatG_{\aff}$ will have powerful applications to automorphic forms on $G$.

Of course, at the moment the above dream remains only a dream; however, in the recent years some interesting results about
representation theory of $G_{\aff}$ over either local or global field have appeared. The purpose of this paper is to survey
some of those results; more precisely, we are going to concentrate on two aspects: study of some particular
Hecke algebras in the local case and the study of Eisenstein series in the global case.
All the results that we are going to discuss generalize well-known results for the group $G$ itself;
however, the generalizations are not always straightforward and some new features appear in the affine case.
\subsection{Hecke algebras and Satake isomorphism}
First let us mention that some version of the above representation theory for $\hatG_{\aff}(\calK)$  was developed in
\cite{Kap-H} and \cite{GK1}-\cite{GK3}. This theory looks promising, but we are not going to discuss it in this paper;
on the other hand, in \cite{BK} we
generalized the Satake isomorphism to the case of $G_{\aff}$. Let us first recall the usual Satake isomorphism
(which can be thought of as the starting point for
Langlands duality mentioned above).

Given an open compact subgroup $K$ in $G(\calK)$ one can consider the Hecke algebra $H(G,K)$ of compactly supported $K$-bi-invariant
distributions on $G(\calK)$ (this is an algebra with respect to convolution). We say that  $H(G,K)$ is a {\it Hecke algebra}
of $G$ with respect to $K$. $H(G,K)$  is a unital  associative algebra and it is
well-known that the study of the  representation theory of $G(\calK)$  is essentially equivalent to studying representation theory
of $H(G,K)$ for different choices of $K$.

The  group
$G(\calO)$ is a maximal compact subgroup of $G(\calK)$. We denote  the corresponding Hecke algebra $H(G,G(\calO))$ by
$\calH_{\sph}(G,\calK)$ and call it the {\em spherical Hecke algebra}. The Satake isomorphism is an
isomorphism between $\calH_{\sph}(G,\calK)$ and the complexified
Grothendieck ring $K_0(\Rep(G^{\vee}))$
of finite-dimensional representations of $G^{\vee}$. For future purposes it will be convenient to note that $K_0(\Rep(G^{\vee}))$
is also naturally isomorphic to the algebra $\CC(T^{\vee})^W$ of polynomial functions on the maximal torus
$T^{\vee}\subset G^{\vee}$ invariant under $W$.

The Satake isomorphism was generalized to $G_{\aff}$ in \cite{BK} and it was studied in more detail in
\cite{BKP}. The formulation is similar to the case of $G$ but the details are somewhat more involved. We shall give
the precise formulation in Section \ref{aff-local}, where we also discuss the analog of the
so called Iwahori-Hecke algebra for $G_{\aff}$ and the affine version of the Gindikin-Karpelevich formula and the
Macdonald formula for the
spherical function (following the papers \cite{BFK}, \cite{BGKP}, \cite{BKP}).

A (partial) analog of the geometric Satake isomorphism (cf. \cite{MV}) in the affine case (following the papers
\cite{BF1}-\cite{BF3}) is also discussed in Section \ref{aff-local}.

\subsection{Eisenstein series for $\hatG_{\aff}$}
The most basic example of automorphic forms, the Eisenstein series for $\hatG_{\aff}$ coming from the Borel subgroup
of $\hatG_{\aff}$, were studied extensively in the works of Garland
(cf. e.g. \cite{Gar1}-\cite{Gar3} and references therein) and Kapranov \cite{Kap-E}. In the forthcoming publication
\cite{BK-Eis} we are going to continue the above study of Eisenstein series
and give some applications (for example, we are going to describe an affine version of the Tamagawa number formula for $G_{\aff}$).
We hope that the above results should have interesting applications to automorphic $L$-functions for the group $G$ itself,
using the affine version of the so called Langlands-Shahidi method (cf. \cite{Gar3}). We discuss it in more detail in
Section \ref{aff-global}.

\subsection{Contents of the paper}
In Section \ref{dig-local} we are going to review some facts about spherical and Iwahori
Hecke algebras of reductive groups over a local non-archimedian field that we are going to generalize
later to the affine case. In addition we also review the corresponding geometric Satake isomorphism.
In Section \ref{dig-global} we
review some well-known facts about automorphic forms and, in particular, about Eisenstein series over a global field $F$.
To simplify the discussion we
concentrate on the case when $F$ is a functional field and the automorphic forms in question are everywhere unramified.

In Section \ref{aff-local} we generalize the constructions of Section \ref{dig-local} to the case of $G_{\aff}$.
Similarly, in Section \ref{aff-global} we generalize the construtions of Section \ref{dig-global} to the affine case.
In particular, we discuss the foundations of the
theory of Eisenstein series for $\hatG_{\aff}$ coming from various parabolic subgroups
of $\hatG_{\aff}$ and relate this theory to some infinite products of automorphic L-functions
for the finite-dimensional group $G$. We also
discuss an affine version of the Tamagawa number formula.

Finally in Section \ref{further} we discuss some other works related to the above constructions as well as future directions of research.
\section{Spherical and Iwahori Hecke algebras in the finite-dimensional case}\label{dig-local}
In this section we recall some well-known facts about representation theory and Hecke algebras
(and their geometric version) of reductive
groups over a local non-archimedian field $\calK$,
which we are going to generalize to the affine case.

\subsection{Groups over local non-archimedian fields and their representations}
For any open compact subgroup $K$ of $G(\calK)$, we denote by $\calH(G,K)$ the Hecke algebra of $G$ with respect to $K$.
A choice of a Haar measure on $G(\calK)$ provides  an identification of
 $\calH(G,K)$ with the space
of $K$-bi-invariant functions on $G(\calK)$.

For every $K$ as above, there is a natural functor from $\calM(G(\calK))$ to the category of left $\calH(G,K)$-modules, sending
every representation $V$ to the corresponding space $V^K$ of $K$-invariants; thus one can try do understand the category $\calM(G,K)$
by studying the categories of $\calH(G,K)$-modules for different $K$.

There are two
choices of $K$ that will be of special interest to us. The first case is $K=G(\calO)$ (which is a maximal
compact subgroup of $G(\calK)$). Recall that the corresponding Hecke algebra in this case is called
{\em the spherical Hecke algebra} and it is
denoted it by $\calH_{\sph}(G,\calK)$. The second case is the case when $K$ is the Iwahori subgroup $I$
(by definition, this is a subgroup of $G(\calO)$ which is equal to the preimage of a Borel subgroup
in $G(\kk)$ under the natural map $G(\calO)\to G(\kk)$).
Let us remind the description of the corresponding algebras in this case.
\subsection{Satake isomorphism}The spherical Hecke algebra $\calH_{\sph}(G,\calK)$ is commutative.
The Satake isomorphism is a canonical
isomorphism between $\calH_{\sph}(G,\calK)$ and the algebra $\CC[\Lam]^W$. The latter algebra has several other standard interpretations: it is also isomorphic to the complexified
Grothendieck ring $K_0(\Rep(G^{\vee}))$
of finite-dimensional representations of $G^{\vee}$ as well as to the algebra $\CC(T^{\vee})^W$
of polynomial functions on the maximal torus
$T^{\vee}\subset G^{\vee}$ invariant under  $W$ of $G$.

The Satake isomorphism is one of the starting points for the celebrated Langlands conjectures.
We are going to present a generalization of the Satake isomorphism to the case of affine Kac-Moody groups in Section \ref{aff-local}.
\subsection{The Iwahori-Hecke algebra}\label{fin-iwahori}
The algebra $\calH(G,I)$ is known to have the following presentation (usually called
"Bernstein presentation"). It is generated by elements $X_{\lam}$ for $\lam\in \Lam$ and $T_w$ for $w\in W$, subject to the following
relations:

1) $T_wT_{w'}=T_{ww'}$ if $\ell(ww')=\ell(w)+\ell(w')$;

2) $X_{\lam}X_{\mu}=X_{\lam+\mu}$; in other words, $X_{\lam}$'s generate the algebra $\CC(T^{\vee})$ inside $\calH(G,I)$;

3) For any $f\in \CC(T^{\vee})$ and any simple reflection $s\in W$ we have
$$
f T_s-T_s s(f)=(q-1)\frac{f-s(f)}{1-X_{-\alp_s}}
$$
(it is easy to see that the right hand side is again an element of $\CC(T^{\vee})$.

The spherical Hecke algebra is a subalgebra of $\calH(G,I)$; in terms of the above presentation it is equal to
$P\cdot\calH(G,I)\cdot P$ where $P=\sum_{w\in W} T_w$.
\subsection{Explicit description of the Satake isomorphism}
Recall, that the Satake isomorphism is an isomorphism $\calS$ between the algebra $\calH_{\sph}(G,\calK)$ and
$\CC(T^{\vee})^W$. Both algebras possess a natural basis, parameterized by the set $\Lam_+$ of dominant coweights of
$G$.

Namely, let $\varpi\in\calO$ be a uniformizer;
for $\lam\in \Lam$ let us denote by $\varpi^{\lam}$ the image of $\varpi\in\calK^*$
under the map $\lam:\calK^*\to T(\calK)\subset G(\calK)$. Then it is known that $G(\calK)$ is the disjoint union of the
cosets $G(\calO)\cdot \varpi^{\lam}\cdot G(\calO)$ when $\lam$ runs over $\Lam_+$.
For every $\lam\in\Lam_+$ we denote by $h_{\lam}\in\calH_{\sph}(G,\calK)$ the characteristic function of the corresponding double coset.
This is a basis of $\calH_{\sph}(G,\calK)$.

On the other hand, for $\lam\in\Lam_+$ let $L(\lam)$ denote the irreducible representation of $G^{\vee}$ with highest weight $\lam$.
Then the characters $\chi(L(\lam))$ form a basis of $\CC(T^{\vee})^W$. We would like to recall what happens to these bases
under the Satake isomorphism (and its inverse).

Let
$W_{\lam}$ is the stabilizer of $\lam$ in $W$ and set
$$
W_{\lam}(q)=\sum\limits_{w\in W_{\lam}} q^{\ell(w)}.
$$

\begin{theorem}[Macdonald, \cite{mac}]\label{macdonald-finite}
For any $\lam\in \Lam^+$ we have
\begin{equation}\label{mac}
\calS(h_{\lam})=\frac{q^{\la\lam,\rho^{\vee}\ra}}{W_{\lam}(q^{-1})}\sum\limits_{w\in W}
w\Bigg(e^{\lam}\frac{\prod\limits_{\alp\in R_+} 1-q^{-1}e^{-\alp}}{\prod\limits_{\alp\in R_+} 1-e^{-\alp}}\Bigg).
\end{equation}
Here $\rho^{\vee}$ is the half-sum of the positive roots of $G$.
\end{theorem}
Let us stress the following: there are many proofs of the above theorem, but essentially all of them use the Iwahori
Hecke algebra $\calH(G,I)$.

On the other hand, we can consider the Satake isomorphism $\calS$ as  is an
isomorphism between the algebras $\calH_{\sph}(G,\calK)$ and
$K_0(\Rep(G^{\vee}))$. For any $\lam\in\Lam_+$ let $L(\lam)$ denote the irreducible representation of
$G^{\vee}$ with highest weight $\lam$.
Then their classes $[L(\lam)]$ form a basis of $K_0(\Rep(G^{\vee}))$. The functions $\calS^{-1}([L(\lam)])$ were described
by Lusztig \cite{Lu-qan} (cf. also \cite{Bry} and \cite{Kato}).
Namely, we have:
\begin{theorem}[Lusztig, \cite{Lu-qan}]\label{qan}
Let $\lam,\mu\in\Lam_+$. Then
$\calS^{-1}([L(\lam)])(\varpi^{\mu})$ is non-zero if and only if $\mu$ is a weight of $L(\lam)$ and in that
case $\calS^{-1}([L(\lam)])(\pi^{\mu})$ is a certain $q$-analog of the weight multiplicity $\dim L(\lam)_{\mu}$
(i.e. it is a polynomial in $q$ with integral coefficients whose value at $q=1$ is equal to the weight
multiplicity).
\end{theorem}
\subsection{Gindikin-Karpelevich formula}\label{gindikin}
The classical Gindikin-Karpelevich formula describes explicitly how a certain intertwining
operator acts on the spherical vector in a principal series representation of $G(\calK)$.
\footnote{More precisely, the Gindikin-Karpelevich
 formula answers the analogous question for real groups} In more explicit terms
it can be formulated as follows.

Let us choose a Borel subgroup $B$ of $G$ and an opposite Borel subgroup
$B_-$; let $U,U_-$ be their unipotent radicals. In addition, let  $L$ denote the coroot lattice of $G$, $R_+\subset L$ -- the set
of positive coroots, $L_+$ -- the subsemigroup of $L$ generated by $R_+$.
Thus any $\gam\in\Lam_+$ can be written as $\sum a_i\alp_i$ where $\alp_i$ are the simple roots.
We shall denote by $|\gam|$ the sum of all the $a_i$.

Set now $\Gr_G=G(\calK)/G(\calO)$. Then it is known that $\calU(\calK)$-orbits on $\Gr$ are in one-to-one correspondence
with elements of $\Lam$;
for any $\mu\in\Lam$ we shall denote by $S^{\mu}$ the corresponding orbit.
The same thing is true for $U_-(\calK)$-orbits. For each $\gam\in\Lam$ we shall denote by
$T^{\gam}$ the corresponding orbit. It is well-known that
$T^{\gam}\cap S^{\mu}$ is non-empty iff $\mu-\gam\in L_+$ and in that case the above intersection
is finite. The Gindikin-Karpelevich formula allows one to compute the number of points
in $T^{-\gam}\cap S^0$ for $\gam\in L_+$ (it is easy to see that the above intersection
is naturally isomorphic to $T^{-\gam+\mu}\cap S^{\mu}$ for any $\mu\in\Lam$). The answer is most easily stated
in terms of the corresponding generating function:
\begin{theorem}[Gindikin-Karpelevich formula]\label{gk-finite}
$$
\sum\limits_{\gam\in\L_+} \# (T^{-\gam}\cap S^0)q^{-|\gam|}e^{-\gam}=
\prod\limits_{\alp\in R_+}\frac{1-q^{-1}e^{-\alp}}{1-e^{-\alp}}.
$$
\end{theorem}
The proof of the above formula is not difficult -- it can be reduced to $G=SL(2)$. However, it can also be obtained
from the Macdonald formula (\ref{mac}) by a certain limiting procedure. The second proof is important
for our purposes since we are going
to use similar argument in the affine case.
\subsection{Geometric Satake isomorphism}\label{geom-satake}
The Satake isomorphism recalled above has a geometric version, called {\em the geometric Satake isomorphism}.
\footnote{As was mentioned above the usual Satake isomorphism is the starting point for Langlands duality;
the classical Langlands duality has its geometric counterpart, usually referred to as
the {\em geometric Langlands duality} which is based on the  geometric version of the Satake isomorphism.}
Let us recall some facts about the geometric Satake isomorphism; later, we are going to discuss a (partial) generalization
of it to the affine case. It is probably worthwhile to note that there exists a geometric approach to the Iwahori-Hecke algebra
(it was developed in the works of Ginzburg, Kazhdan-Lusztig and Bezrukavnikov), but we are not going to discuss it in this paper.

Let now $\calK=\CC((s))$ and let $\calO=\CC[[s]]$;
here $s$ is a formal variable. Let $\Gr_G=G(\calK)/G(\calO)$.
Then the geometric (or categorical)
analog of the algebra $\calH(G)$ considered above is the category $\Perv_{G(\calO)}(\Gr_G)$ of
$G(\calO)$-equivariant perverse sheaves on $\Gr_G$ (cf. e.g. [27]). According to {\em loc. cit.} the category $\Perv_{G(\calO)}(\Gr_G)$
possesses canonical tensor structure and the geometric Satake isomoprhism asserts that this category is
equivalent to $\Rep(G^{\vee})$ as a tensor category. The corresponding fiber functor from $\Perv_{G(\calO)}(\Gr_G)$
to vector spaces sends every perverse sheaf $\calS\in\Perv_{G(\calO)}(\Gr_G)$ to its cohomology. Another way to construct
this ("fiber") functor is discussed in [32].

More precisely, one can show that $G(\calO)$-orbits on $\Gr_G$ are finite-dimensional
and they are indexed by the set $\Lam^+$ of dominant weights of $G^{\vee}$. For every $\lam\in\Lam^+$ we denote
by $\Gr_G^{\lam}$ the corresponding orbit and by $\oGr_G^{\lam}$ its closure in $\Gr_G$. Then $\Gr_G^{\lam}$ is a
non-singular quasi-projective algebraic variety over $\CC$ and $\oGr_G^{\lam}$ is a (usually singular)
projective variety. One has
$$
\oGr_G^{\lam}=\bigcup\limits_{\mu\leq \lam}\Gr_G^{\lam}.
$$
One of the main properties of the geometric Satake isomorphism is that it sends the irreducible $G^{\vee}$-module
$L(\lam)$ to the intersection cohomology complex $\IC(\oGr_G^{\lam})$. In particular, the module $L(\lam)$
itself gets realized as the intersection cohomology of the variety $\oGr_G^{\lam}$.

As a byproduct of the geometric Satake isomorphism one can compute $\IC(\oGr_G^{\lam})$ in terms of the
module $L(\lam)$. Namely, it is well-known that the stalk
of $\IC(\oGr_G^{\lam})$ at a point of $\Gr_G^{\mu}$ as a graded vector space is essentially
equal to the associated graded $\gr^F L(\lam)_{\mu}$ of the $\mu$-weight space $L(\lam)_{\mu}$ in $L(\lam)$
with respect to certain filtration, called the Brylinski-Kostant filtration. This is a geometric analog of Theorem
\ref{qan}.

One can construct certain canonical transversal slice $\ocalW^{\lam}_{G,\mu}$ to
$\Gr_G^{\mu}$ inside $\oGr_G^{\lam}$. This is a conical affine algebraic variety (i.e. it is endowed with an
action of the multiplicative group $\GG_m$ which contracts it to one point).
The above result about the stalks of $\IC(\oGr_G^{\lam})$
then gets translated into saying that the stalk of the IC-sheaf of $\ocalW^{\lam}_{G,\mu}$ at the unique $\GG_m$-fixed point
is essentially isomorphic to $\gr^F L(\lam)_{\mu}$. Note that since  $\ocalW^{\lam}_{G,\mu}$ is contracted to the above point by the
$\GG_m$-action, it follows that the stalk of of $\IC(\ocalW^{\lam}_{G,\mu})$ is equal to the global intersection cohomology
$\IH^*(\ocalW^{\lam}_{G,\mu})$.

The varieties $\ocalW^{\lam}_{G,\mu}$ are important for us because we are going to describe their analogs when
$G$ is replaced by $G_{\aff}$; on the other hand, the affine analogs of
$\oGr^{\lam}_G$ should be wildly infinite-dimensional and
we don't know how to think about them.
\section{Unramified automorphic forms in the finite-dimensional case}\label{dig-global}
In this Section we recall some classical facts about automorphic forms on reductive groups
of reductive groups over a global field, which we are going to generalize to the affine case later.  To simplify the discussion,
we  restrict ourselves to unramified automorphic forms over a functional field $F$ (i.e. the field of rational functions
on a smooth projective curve $X$ over a finite field $\kk$). In this case, there is a way to think
about automorphic forms in terms of functions on the moduli space of $G$-bundles over the curve $X$; this
point of view  will be convenient for us in the affine case.
\footnote{We would like to note that the pioneering work on the subject in the  affine case was done by H.~Garland
\cite{Gar1}-\cite{Gar3}, who dealt with the case when $F=\QQ$.
Most of the results that we are going to describe in the affine case can be generalized
to any global field, but technically the case of number fields is more difficult
due to the existence of archimedian places and we
prefer not to discuss it in this survey paper}
\subsection{Automorphic forms over functional fields and $G$-bundles}\label{aut-finite}
Let $X$ be a smooth projective geometrically irreducible curve over
a finite filed $\kk=\FF_q$. Let also $G$ be a split semi-simple simply connected group over $\kk$.
We set $F=\kk(X)$; this is a global field and we let $\AA_F$ denote its ring of adeles. We also denote by $\calO(\AA_F)$ ring of integral adeles.

It is well-know that the double quotient
\begin{equation}\label{bun=adeles}
G(\calO(\AA_F))\backslash G(\AA_F)/G(F)\simeq \Bun_G(X).
\end{equation}
So far, this is just an isomorphism between two sets, but later we are going to upgrade it to an an equivalence between two groupoids.
(Complex valued) functions on the above space will usually be referred to
as unramified automorphic forms. We denote the space of such functions by $\CC(\Bun_G(X))$.
\subsection{The Hecke algebra action}Recall, that for a local field $\calK$ with ring of integers $\calO_{\calK}$ we denote by $\calH_{\calK}(G)$
the {\em spherical Hecke algebra} of $G$ over $\calK$; this algebra is isomorphic
to the algebra $\CC[G^{\vee}]^{G^{\vee}}$ of ad-invariant regular functions on $G^{\vee}$.

Let $v$ be a place of $F$ and let $\calK_v$ be the corresponding local completion of $F$ with ring of integers $\calO_v$.
Let $q_v$ be the number of elements in the residue field of $\calO_v$. Let $\calH_v(G)$ denote the corresponding Hecke algebra.
It is well-known that $\calH_v(G)$ acts on $\CC(\Bun_G(X))$ by correspondences.

Assume that $f\in \CC(\Bun_G(X))$ is an eigen-function
of all the $\calH_v(G)$. We say that $f$ is a Hecke eigen-form. The corresponding eigen-value is an element
$g_v(f)\in G^{\vee}$ for all places $v$. Given a finite-dimensional
representation $\rho:G^{\vee}\to \Aut(V)$, we define the $L$-function $L(f,\rho,s)$ of $f$
by
$$
L_G(f,\rho,s)=\prod\limits_{v} \frac{1}{\det(1- q_v^{-s}\rho(g_v))}.
$$
\subsection{Groupoids}In what follows it will be convenient to treat various double coset spaces
as groupoids rather than as sets.
By a groupoid we shall mean a small category $\calX$ where all morphisms are isomorphisms.
Any such groupoid is equivalent to a quotient groupoid groupoid $X/H$ where $X$ is a set and $H$ is a group acting on $X$.
Given a group $G$ and two subgroups $H_1$ and $H_2$, the double quotient
$H_1\backslash G/H_2$ has a natural structure of a groupoid; in particular, this applies to the double coset space
(\ref{bun=adeles}).

We shall sometimes denote by $|\calX|$ the set of isomorphism classes of objects of $\calX$. Given $x\in\calX$
we denote by $\Aut_{\calX}(x)$ the automorphism group of $x$ in $\calX$. Given two groupoids $\calX$ and $\calY$ it makes
sense to consider functors $f:\calX\to\calY$. For every $y\in Y$,  the fiber $f^{-1}(y)$ is also a groupoid:
by definition
$\Aut_{f^{-1}(y)}(x)$ is the group of automorphisms of $x$ whose image in the group of automorphisms of
$y$ is trivial. We say that $f$ is representable if $\Aut_{f^{-1}(y)}(x)$ is trivial for every $x$ and $y$.

For a map $p:X\to Y$ between two sets we are going to denote by $p^*$ the pull-back of functions (this is a linear map from
functions on $Y$ to functions on $X$) and by $p_!$ the operation of summation over the fiber (this is a linear map
from functions on $X$ to functions on $Y$; a priori it is well-defined when $p$ has finite fibers). More generally, we
can talk about $p^*$ and $p_!$ when $X$ and $Y$ are groupoids. In this case we should define $p_!$ in the following
way
\begin{equation}\label{pshrieck}
p_!(f)(y)=\sum\limits_{x\in |p^{-1}(y)|} \frac{f(x)}{\# \Aut_{p^{-1}(y)}(x)}.
\end{equation}
In particular, we can apply this to the case when $p$ is the map from $\calX$ to the point. In this case, we shall
denote $p_!$ by $\int_{\calX}$. In other words, for a function $f:|\calX|\to \CC$ we set
\begin{equation}\label{measure}
\int\limits_{\calX} f=\sum\limits_{x\in |\calX|} \frac{f(x)}{\# \Aut (x)}.
\end{equation}
The above sum makes sense if both $|\calX|$ and $\Aut(x)$ are finite. When $|\calX|$ is infinite, the sum sometimes still makes
sense (if it converges).

\subsection{Eisenstein series and constant term}\label{e-finite}Let $P\subset G$ be parabolic subgroup and let $M$ be the corresponding Levi group.
We have canonical maps $G\leftarrow P\to M$ which give rise to the diagram
$$
\begin{CD}
\Bun_P&(X)@>\eta>>\Bun_M(X)\\
@V\pi VV\\
\Bun_G&(X)
\end{CD}
$$

The map $\pi$ is representable but has infinite fibers. The map $\eta$ has finite fibers but it is not representable.

The connected components of $\Bun_M$ are numbered by elements of the lattice $\Lam_M=\Hom(\GG_m,M/[M,M])$; we denote
by $2\rho_P$ the element of the dual lattice $\Hom(M,\GG_m)$, equal to the determinant of $M$-action on $\grn_{P}$.
Abusing the notation, we shall denote by the same symbol the corresponding function
$\Bun_M(X)\to \ZZ$; it also makes sense to talk about the function $\rho_P:\Bun_M(X)\to \ZZ$.
\footnote{$\rho_P$ takes values in $\ZZ$ and not $\frac{1}{2}\ZZ$ since $G$ was assumed to be simply connected}

Given a function $f\in \CC(\Bun_M(X))$ we define the Eisenstein series $\Eis_{GP}(f)\in \CC(\Bun_G(X))$ by setting
$$
\Eis_{GP}(f)=\pi_!\eta^*(fq^{\rho_P}).
$$
This makes sense when $f$ has finite support. When $f$ has infinite support, well-known convergence issues arise.

Similarly, for a function $g\in\CC(\Bun_G(X))$ we define the constant term $c_{GP}(g)\in \CC(\Bun_M(X))$
by
$$
c_{GP}(g)=q^{\rho_P}\eta_! \pi^*(g).
$$
The constant term is well-defined for any $g$. We say that $g$ is cuspidal if $c_{GP}(g)=0$ for all $P$.
\subsection{Constant term of Eisenstein series}Let us recall how to compute the composition of Eisenstein series and constant
term operators. For simplicity we are going to restrict ourselves to the following situation.
Let $P$ and $M$ be as above and let $P_-$ be a parabolic subgroup opposite to $P$ (i.e. $P\cap P_-=M$).
Let $f$ be a cuspidal function on $\Bun_M$.
Let $M^{\vee}$ be the corresponding Langlands dual group and let $P^{\vee}$, $P_-^{\vee}$ be the corresponding
parabolics in $G^{\vee}$; let also $\grn_{P}^{\vee}$, $\grn_{P_-}^{\vee}$ be the nilpotent radicals of their
Lie algebras. The group $M^{\vee}$ acts on both of these space. For each $w\in N(M)/M=N(M^{\vee})/M^{\vee}$
we denote by by $P^w, (P^{\vee})^w, \grn_P^w$ etc. the corresponding groups or Lie algebras obtained. Similarly,
for $f\in \Bun_M(X)$ we denote by $f^w$ the corresponding function on $\Bun_M$ obtained by applying conjugation by
$w$ to $f$.

Then we have the following well known:
\begin{theorem}\label{gk}
Assume that $f\in \CC(\Bun_M(X))$ is a cuspidal Hecke eigen-function.
\begin{enumerate}
\item
$$
c_{GP}\circ \Eis_{GP}(f)=\sum\limits_{w\in N(M)/M} q^{(g-1)\dim \grn_P^{\vee}\cap (\grn_P^{\vee})^w}f^w
\frac{L(f,\grn_P^{\vee}\cap (\grn_{P_-}^{\vee})^w, 0)}{L(f,\grn_P^{\vee}\cap (\grn_{P_-}^{\vee})^w, 1)}.
$$
\item
$$
c_{GP_-}\circ \Eis_{GP}(f)=\sum\limits_{w\in N(M)/M} q^{(g-1)\dim \grn_P^{\vee}\cap (\grn_{P_-}^{\vee})^w}f^w
\frac{L(f,\grn_P^{\vee}\cap (\grn_{P}^{\vee})^w, 0)}{L(f,\grn_P^{\vee}\cap (\grn_{P}^{\vee})^w, 1)}.
$$
\end{enumerate}
\end{theorem}

\noindent
{\bf Remark.} Similar formula holds for the composition $c_{GQ}\circ \Eis_{GP}$ where $Q$ is any parabolic having $M$ as its
Levi subgroup. However, only the above cases will be relevant for us in the sequel.

\subsection{Tamagawa number}
In this section we recall the calculation of the Tamagawa number in the finite-dimensional case.
For simplicity we stick again to the functional field case; we also assume that $G$ is simply connected.

Let $X$, $G$, $\Bun_G(X)$ be as above.
Let $\zeta(s)$ denote the $\zeta$-function of the curve $X$. In other words,
$$
\zeta(s)=\frac{\det(1-q^{-s}\Fr:H^1(X,\qlb)\to H^1(X,\qlb))}{(1-q^{-s})(1-q^{-s+1})}
$$
Let also $d_1,\cdots,d_r$ be the degrees of the generators of the ring of invariant polynomials
on $\grg=Lie(G)$ (also called the exponents of $G$).
\begin{theorem}\label{tam-finite}
We have
\begin{equation}\label{tamagawa-finite}
\int\limits_{\Bun_G(X)} 1=q^{(g-1)\dim G}\prod\limits_{i=1}^l \zeta(d_i).
\end{equation}
\end{theorem}
The above formula was proved by Langlands in \cite{Lan-tam}\footnote{Langlands considered the case of global field $\QQ$ but
the argument in \cite{Lan-tam} applies to any global field}. It can also be derived from a computation
of $H^*(\Bun_G,\qlb)$ (due to a result of Behrend \cite{Beh}),
which was done  in \cite{HS} (following a computation of Atiyah-Bott \cite{AB} over $\CC$). Let us formulate this result,
since this interpretation of Theorem \ref{tam-finite} is instructive from the point of view of generalization to the affine case.

It is well-known that $H^*(pt/G,\qlb)=\Sym^*(V)$ for a certain graded vector space endowed with an action of Frobenius.
Explicitly, $V$ has a basis $c_1,...,c_{r}$ where each $c_i$ has cohomological degree $2d_i$ and Frobenius eigenvalue
$q^{d_i}$.
\begin{theorem}\label{at-bott}
$H^*(\Bun_G,\qlb)=\Sym^*(V\ten H_*(X,\qlb))$; here $\Sym^*$ is understood in the "super"-sense and $H_*(X,\qlb)$ stands for
the homology of $X$ defined by setting $H_i(X,\qlb)$ to be the dual space to $H^{-i}(X,\qlb)$ (thus $H_*$ "sits" in degrees
$0,-1,-2$).
\end{theorem}

\subsection{Sketch of the proof of Theorem \ref{tam-finite}}\label{sketch}
Let us recall the basic steps of Langlands' proof of Theorem \ref{tam-finite}, since this proof will actually be important for
some definitions in the affine case.

For any $\gam\in \Lam$ let us denote by $p^{\gam}$ the restriction of the natural map
$p:\Bun_B\to \Bun_G$ to $\Bun_B^{\gam}$. Let us denote by $\grt^*$ the dual space to the Lie algebra of $T$ (over $\CC$).
Then for any $\grs\in \grt^*$ let us set
$$
E_{\gam}=(p_{\gam})_!({\mathbf 1}_{\Bun_B^{\gam}});
\quad E(\grs)=\Eis_{GB}(\sum \limits_{\gam}{\mathbf 1}_{\Bun_T^{\gam}}q^{\la\gam,\grs\ra})=
\sum\limits_{\gam\in\Lam}q^{\la\gam,\grs+\rho\ra}E_{\gam}.
$$
The following result is well-known:
\begin{lemma}
\begin{enumerate}
\item
The above series is absolutely convergent for $\Re(\grs)>\rho^{\vee}$ and it extends to a meromorphic function on the whole
of $\grt^*$.
\item
The function $E(\grs)$ has simple pole along every hyperplane of the form $\la \grs-\rho,\alp^{\vee}\ra=0$ where $\alp$ is a
simple coroot and no other poles near $\grs=\rho$.
\item
The residue of $E(\grs)$ at $\grs=\rho$ defined as the limit
$
\lim\limits_{\grs\to \rho}E(\grs) \prod\limits_{i\in I} \la \grs-\rho,\alp_i^{\vee}\ra
$
is a constant function on $\Bun_G$.
\end{enumerate}
\end{lemma}
Let us call the above residue $\grr(G)$.
The idea of the proof of Theorem \ref{tam-finite} is this: on the one hand it is easy to compute
$\int\limits_{\Bun_G}\grr(G)\cdot {\mathbf 1}_{\Bun_G}$.
On the other hand, one can compute $\grr(G)$ itself by computing the constant term of $\grr(G)\cdot {\mathbf 1}_{\Bun_G}$.
More precisely, one can deduce Theorem \ref{tam-finite} from the following
\begin{theorem}\label{finite-residue}
\begin{enumerate}
\item
$\int\limits_{\Bun_G} \grr(G)\cdot {\mathbf 1}_{\Bun_G}=q^{(g-1)\dim N}(\ln q)^{-r} ~\frac{(\# \Pic^0(X))^{r}}{(q-1)^{r}}.$
\item
$c_{GB}(\grr(G)\cdot {\mathbf 1}_{\Bun_G})|_{\Bun_T^{\gam}}=q^{-\la\gam,\rho\ra}~ \frac{(\Res_{s=1}\zeta(s))^{r}}
{\prod\limits_{i=1}^{r} \zeta (d_i)}.$
\end{enumerate}
\end{theorem}
Theorem \ref{finite-residue} implies Theorem \ref{tam-finite} since
$
\Res\limits_{s=1}\zeta(s)=\frac{\# \Pic^0(X)}{q^{g-1}(q-1)\ln q}.
$

\section{Spherical and Iwahori Hecke algebras of affine Kac-Moody groups over a local non-archimedian field}\label{aff-local}

In this Section we discuss analogs of spherical and Iwahori Hecke algebra for the group $G_{\aff}$ over a local non-archimedian field. The results are taken from \cite{BK}, \cite{BFK}, \cite{BGKP} and \cite{BKP}.
\subsection{The semi-group $G^+_{\aff}(\calK)$ and the affine spherical Hecke algebra}
One may consider the group $G_{\aff}(\calK)$ and its subgroup
$G_{\aff}(\calO)$.
The group $G_{\aff}$ by definition maps to $\GG_m$; thus $G_{\aff}(\calK)$
maps to $\calK^*$. We denote by $\varrho$ the composition of this map with the valuation maps$\cal^*\to \ZZ$.

We now define the semigroup $G_{\aff}^+(\calK)$ to be the subsemigroup of $G_{\aff}(\calK)$ denerated by:

$\bullet$ the central $\calK^*\subset G_{\aff}(\calK)$;

$\bullet$ the subgroup $G_{\aff}(\calO)$;

$\bullet$ All elements $g\in G_{\aff}(\calK)$ such that $\varrho(g)>0$.

We show in \cite{BK} that there exists an associative algebra structure on a  completion $\calH_{\sph}(G_{\aff},\calK)$  of the space of finite linear combinations of
double cosets of $G_{\aff}(\calK)^+$ with respect to $G_{\aff}(\calO)$. We would like to emphasize that this statement is by no means trivial -- in \cite{BK} it is proved using some cumbersome algebro-geometric machinery.
A more elementary proof of this fact is going to appear in \cite{BKP} but that proof is still rather long.

We call the above algebra {\it the spherical Hecke algebra of $G_{\aff}$}.
The algebra $\calH_{\sph}(G_{\aff},\calK)$ is graded by non-negative integers
(the grading comes from the map $\varrho$ which is well-defined on double cosets
with respect to $G_{\aff}(\calO)$); it is also an algebra over the field $\CC((v))$ of Laurent power series in a variable
$v$, which comes from the central $\calK^*$  in $G_{\aff}(\calK)$.
\subsection{The affine Satake isomorphism}
The statement of the Satake isomorphism for $G_{\aff}$ is very similar to that for $G$.
First of all, in \cite{BK} we define an analog of the algebra $\CC(T^{\vee})^W$ which we shall
denote by $\CC(\hatT^{\vee}_{\aff})^{W_{\aff}}$ (here $T_{\aff}=\CC^*\x T^{\vee}\x \CC^*$ is the dual of the maximal
torus of $G_{\aff}$, $W_{\aff}$ is the corresponding affine Weyl group and $\CC(\hatT^{\vee})$ denotes certain completion
of the algebra of regular functions on $T^{\vee}_{\aff}$). This is a finitely generated $\ZZ_{\geq 0}$-graded commutative algebra
over the field $\CC((v))$ of Laurent formal power series in the variable $v$ which should be thought of as a coordinate on
the third factor in $T^{\vee}_{\aff}=\CC^*\x T^{\vee}\x\CC^*$ (the grading has to do with the first factor);
moreover, each component of the
grading is finite-dimensional over $\CC((v))$.

To simplify notations we will always assume that $G$ is simple and  simply connected (although the case when $G$ is a torus
is also very instructive and it is much less trivial than in the usual case - cf. Section 3 of \cite{BK}).

In this case we define (in \cite{BK}) the {\it Langlands dual group} $G_{\aff}^{\vee}$, which is a group ind-scheme over $\CC$.
Then
$G_{\aff}^{\vee}$ is another Kac-Moody group whose Lie algebra $\grg_{\aff}^{\vee}$ is an affine
Kac-Moody algebra with root system dual to that of $\grg_{\aff}$ (thus, in particular, it might be a twisted affine Lie algebra when
$\grg$ is not simply laced).
The group $G_{\aff}^{\vee}$ contains the torus $T_{\aff}^{\vee}$; moreover the first $\CC^*$-factor in $T^{\vee}_{\aff}$ is central
in $G^{\vee}_{\aff}$; also the projection $T^{\vee}_{\aff}\to\CC^*$ to the last factor extends to a homomorphism
$G_{\aff}^{\vee}\to \CC^*$.
One defines a category $\Rep(G_{\aff}^{\vee})$ of
$G_{\aff}^{\vee}$-modules which properly contains all highest weight integrable representations of finite length and
also certain infinite direct sums of irreducible highest weight integrable  representations which
which is stable under tensor product. The  character map provides
an isomorphism of the complexified Grothendieck ring $K_0(G_{\aff}^{\vee})$ with
the algebra $\CC(\hatT^{\vee}_{\aff})^{W_{\aff}}$. The corresponding grading
on $K_0(G_{\aff})$ comes from the central charge of $G_{\aff}^{\vee}$-modules and the action of the variable $v$
comes from tensoring $G_{\aff}^{\vee}$-modules by the one-dimensional representation coming from the homomorphism
$G_{\aff}^{\vee}\to\CC^*$, mentioned above.

The affine Satake isomorphism (proved in \cite{BK} and reproved in a different way in \cite{BKP}) claims the following:
\begin{theorem}\label{aff-satake}
The algebra  $\calH_{\sph}(G_{\aff},\calK)$ is canonically isomorphic
to $\CC(\hatT^{\vee}_{\aff})^{W_{\aff}}$ (and thus also to $K_0(G_{\aff}^{\vee})$).
\end{theorem}

As was mentioned, in the case when $G$ is semi-simple and simply connected, the group $G_{\aff}$ is an affine Kac-Moody group.
We expect that with
slight modifications our Satake isomorphism should make sense for any symmetrizable Kac-Moody group. However, our proofs
are really designed for the affine case and do not seem to generalize to more general Kac-Moody groups. The corresponding
Hecke algebra has recently been defined in \cite{GaRo}.
\subsection{The Iwahori-Hecke algebra of $G_{\aff}$}\label{aff-iwahori}
Let now $I_{\aff}\subset G_{\aff}(\calO)$ be the  Iwahori subgroup and
let $\calH(G_{\aff}^+,I_{\aff})$ denote the space of $I_{\aff}$-bi-invariant functions on $G^+_{\aff}(\calK)$
supported on a union of finitely many double cosets.
It is shown \cite{BKP}  that the usual convolution is well-defined on $\calH(G_{\aff}^+,I_{\aff})$. Note that
this is different from the spherical case, the convolution was only defined on a completion of the space of functions on double cosets with finite support.

The structure of the algebra $\calH(G_{\aff}^+,I_{\aff})$ is similar to the finite-dimensional case.
For simplicity let us assume that $G$ is semi-simple and simply connected. Let
$\Lam_{\aff}=\ZZ\oplus\Lam\oplus\ZZ$ be the lattice of cocharacters of $T_{\aff}$. Let $\Lam_{\aff}^+\subset \Lam$
be the Tits cone, consisting of all elements $(a,\olam,k)\in \Lam_{\aff}$ such that either $k>0$ or
$k=0$ and $\lam=0$. Let $\bfH_{\aff}$ denote the algebra generated by elements $X_{\lam}$, $\lam\in \Lam_{\aff}$ and $T_w$,
$w\in W_{\aff}$ with relations 1), 2), 3) as in Subsection \ref{fin-iwahori} with $q$ replaced by a formal variable $v$. The algebra $\bfH_{\aff}$ is $\ZZ$-graded; this grading is defined by setting
$$
\deg T_w=0;\quad \deg X_{(a,\lam,k)}=k.
$$
We denote by $\bfH_{\aff,k}$ the space of all elements of degree $k$ in $\bfH_{\aff}$. Note that $\bfH_{\aff,0}$ is a subalgebra
of $\bfH_{\aff}$, which is isomorphic to Cherednik's double affine Hecke algebra.

On the other hand, let
$$
\bfH_{\aff}^+=\CC\<T_w\>_{w\in W_{\aff}}\oplus\Bigg(\bigoplus\limits_{k>0}\bfH_{\aff,k}\Bigg).
$$
(it is easy to see that $\bfH_{\aff}^+$ is just generated
by elements $X_{\lam}$, $\lam\in \Lam_{\aff}^+$ and $T_w$, $w\in W$.

The following result is proved in \cite{BKP}:
\begin{theorem}\label{aff-iw-theorem}
The algebra $\calH(G_{\aff}^+,I_{\aff})$ is isomorphic to the specialization of $\bfH_{\aff}^+$ to $v=q$.
\end{theorem}
In particular, the algebra $\calH(G_{\aff}^+,I_{\aff})$ is closely related to Cherednik's double affine Hecke algebra.
We would like to note that another relation between the double affine Hecke algebra and the group $G_{\aff}(\calK)$ was
studied by Kapranov (cf. \cite{Kap-H}).

By definition the algebra $\calH(G_{\aff}^+,I_{\aff})$ by definition is endowed with a natural basis corresponding to characteristic
functions of double cosets of $I_{\aff}$ on $G_{\aff}^+(\calK)$. It is natural to conjecture that this basis comes in fact
from a $\CC[v,v^{-1}]$-basis in $\bfH_{\aff}^+$ but we don't know how to prove this. It would be interesting to give an algebraic
description of this basis.
\subsection{Explicit description of the affine Satake isomorphism}
We would like now to describe the affine analog of Theorem \ref{macdonald-finite} (the proof is going to appear in \cite{BKP}
and it uses the algebra $\calH(G_{\aff}^+,I_{\aff})$ in an essential way). For simplicity we are going to assume again that
$G$ is semi-simple and simply connected.

The (topological) basis $\{ h_{\lam}\}$ of the algebra $\calH_{\sph}(G_{\aff},\calK)$ is defined exactly as in the
case of $G$ and one might expect that
$$
\calS_{\aff}(h_{\lam})=\frac{q^{\la\lam,\rho_{\aff}^{\vee}\ra}}{W_{\aff,\lam}(q^{-1})}\sum\limits_{w\in W_{\aff}}
w\Bigg(e^{\lam}\frac{\prod\limits_{\alp\in R_{+,\aff}} 1-q^{-1}e^{-\alp}}{\prod\limits_{\alp\in R_{+,\aff}} 1-e^{-\alp}}\Bigg)^{m_\alp}.
$$
Here $R_{+,\aff}$ denote the set of positive coroots of $G_{\aff}$, $m_\alp$ is the multiplicity of the coroot $\alp$ and
$\rho^{\vee}_{\aff}$ is the corresponding affine analog of $\rho^{\vee}$.
However, it turns out that the above formula is wrong! Let us explain how to see this (and then we are going to present the correct statement).

Let $\calS_{\aff}:\calH_{\sph}(G_{\aff},\calK)\to \CC(\hatT^{\vee}_{\aff})^{W_{\aff}}$ denote the affine Satake isomorphism.
We would like to compute $\calS_{\aff}(h_{\lam})$ for every $\lam\in\Lam_+$. Let us consider the case $\lam=0$. In this case the element $h_0$ is the unit element of the algebra $\calH_{\sph}(G_{\aff},\calK)$
and thus we must have $\calS_{\aff}(h_0)=1$. On the other hand, the right hand side of (\ref{mac}) in Theorem \ref{macdonald-finite}
in this case equal to
$$
\frac{1}{W_{\aff}(q^{-1})}\sum\limits_{w\in W_{\aff}}
w\Bigg(\frac{\prod\limits_{\alp\in R_{\aff,+}} 1-q^{-1}e^{-\alp}}{\prod\limits_{\alp\in R_{\aff,+}} 1-e^{-\alp}}\Bigg)^{m_{\alp}}.
$$
Let us denote this infinite sum by $\Del$.
This function was studied by Macdonald in \cite{mac-new}. It is easy to see that in the finite case this sum is equal $1$
but in the affine case it is  different from $1$ and  Macdonald gave the following explicit product formula for  $\Del$ which
is going to be important when we discuss affine Eisenstein series. For simplicity let us
assume that $G$ is simply laced and that its Lie algebra is simple.
Let also $\del$ denote the minimal positive imaginary coroot of $G$ and recall that we denote by $d_1,...,d_r$ the degrees
of the generators of the ring $\Sym(\grt)^W$. Then Macdonald proved that
\begin{equation}\label{delta}
\Del=\prod\limits_{i=1}^r\prod\limits_{j=1}^{\infty} \frac{1-q^{-d_i}e^{-j\del}}{1-q^{-d_i+1}e^{-j\del}}.
\end{equation}
Let us now go back to the description of $\calS_{\aff}(h_{\lam})$. The following theorem is proved in \cite{BKP}:
\begin{theorem}\label{aff-macdonald}
For every $\lam\in\Lam_+$ we have
\begin{equation}\label{mac-aff}
\calS_{\aff}(h_{\lam})=\frac{1}{\Del}\cdot\frac{q^{\la\lam,\rho_{\aff}^{\vee}\ra}}{W_{\aff,\lam}(q^{-1})}\sum\limits_{w\in W_{\aff}}
w\Bigg(e^{\lam}\frac{\prod\limits_{\alp\in R_{+,\aff}} 1-q^{-1}e^{-\alp}}{\prod\limits_{\alp\in R_{+,\aff}} 1-e^{-\alp}}\Bigg)^{m_\alp}
\end{equation}
\end{theorem}
The appearance of $\Del$ in the above formula is very curious. Some geometric explanation for it was given in \cite{BFK}
(it is mentioned in a little more detail in the next Subsection).
\subsection{Affine Gindikin-Karpelevich formula}
Let us recall the notations of Subsection \ref{gindikin}; it is clear that at least {\em set theoretically} it makes
sense to consider $\Gr_{\widehat\calG}$ for any completed symmetrizable Kac-Moody group $\widehat\calG$.
The notations $\Lam,\Lam_+,R_+,\Gr_G,S^{\mu},T^{\gam}$
make sense for $\widehat\calG$  without any changes (at least if we think about $S^{\mu}$ and $T^{\gam}$ as sets and not as geometric objects).
\begin{conjecture}\label{gin-fin}
For any $\gam\in\Lam_+$ the intersection $T^{-\gam}\cap S^0$ is finite.
\end{conjecture}
This conjecture is proved in \cite{BFK} when $\calK=\FF_q((t))$ and in \cite{BGKP} for any non-archimedian local field $\calK$ when $\widehat\calG=\hatG_{\aff}$.
\footnote{The general case is probably provable by using the techniques of \cite{GaRo}.}
In the cases when Conjecture \ref{gin-fin} is known, it makes sense to ask whether one can
compute the generating function
$$
I_{\widehat\calG}(q)=\sum\limits_{\gam\in\Lam_+}\# (T^{-\gam}\cap S^0)~q^{-|\gam|}e^{-\gam}.
$$
We do not know the answer in general; however, it \cite{BFK} and \cite{BGKP} it is proved that
when $\widehat\calG=\hatG_{\aff}$ we have
\begin{equation}\label{gin}
I_{\widehat\calG}(q)=\frac{1}{\Del}
\Bigg(\frac{\prod\limits_{\alp\in R_{+,\aff}} 1-q^{-1}e^{-\alp}}{\prod\limits_{\alp\in R_{+,\aff}} 1-e^{-\alp}}\Bigg)^{m_\alp}
\end{equation}
The formula (\ref{gin}) is obtained in \cite{BGKP} as a limit of (\ref{mac-aff}). Earlier, a different proof was given
in \cite{BFK} in the case when $\calK$ is a functional field, using the geometry of Uhlenbeck spaces studied in
\cite{bfg}. We remark that the latter proof also gave a geometric explanation for the appearance of the correction
term $\frac{1}{\Del}$ in this formula. We don't have enough room for details, but
let us just note that  very roughly speaking it is related to the fact that affine Kac-
Moody groups over a functional local field can be studied using various
moduli spaces of bundles on algebraic surfaces.
In particular, by combining the results of \cite{BFK}, \cite{BGKP} and \cite{BKP} one gets a new proof of the identity
(\ref{delta}) (which is independent of Cherednik's results).

\subsection{Towards geometric Satake isomorphism for $G_{\aff}$}\label{aff-geom-satake}
Some parts of the geometric Satake isomorphism have been generalized to $G_{\aff}$ in the papers
\cite{BF1}-\cite{BF3}. The idea the approach of \textit{loc. cit.} belongs in fact to I.~Frenkel who suggested that integrable
representations of $G_{\aff}^{\vee}$ of level $k$ should be realized geometrically in terms of some moduli spaces
related to $G$-bundles on $\AA^2/\Gam_k$, where $\Gam_k$ is the group of roots of unity of order $k$ acting
on $\AA^2$ by $\zeta(x,y)=(\zeta x,\zeta^{-1}y)$.
 \cite{BF1} constitutes an attempt to make this idea more precise.

Let $\Bun_G(\AA^2)$ denote the moduli space of principal $G$-bundles on $\PP^2$ trivialized at the ``infinite"
line
$\PP^1_{\infty}\subset \PP^2$. This is an algebraic variety which has connected components parametrized by non-negative
integers, corresponding to different values of the second Chern class of the corresponding bundles.
Similarly, one can define $\Bun_G(\AA^2/\Gam_k)$. Very vaguely, the main idea
of \cite{BF1} can be formulated in the following way:

\medskip
\noindent
{\bf The basic principle:}

1) The integrable representations of $G_{\aff}^{\vee}$  of level $k$
have to do with the geometry (e.g. intersection cohomology) of some varieties closely related to
$\Bun_G(\AA^2/\Gam_k)$.

2) This relation should be thought of as similar to the relation between finite-dimensional representations
of $G^{\vee}$ and the geometry of the affine Grassmannian $\Gr_G$.

\medskip
\noindent
We believe that 1) above has many different aspects. \cite{BF1} is concentrated on just one such aspect; namely, it is
explained in \cite{BF1}
 how one can construct an analog of the varieties $\ocalW^{\lam}_{\mu}$
in the affine case (using the variety $\Bun_G(\AA^2/\Gam_k)$ as well as the corresponding {\em Uhlenbeck
compactification} of the moduli space of $G$-bundles - cf. \cite{bfg}).
It is conjectured that the stalks of IC-sheaves of these varieties are governed by the affine version of
$\gr^F L(\lam)_{\mu}$.
\footnote{In fact, the definition of the filtration $F$ given in \cite{BF1} is slightly wrong;
the correct definition is given in \cite{slofstra}.}
This conjecture is still open, but in \cite{BF1} it is proved
in several special cases. More precisely, it is proved in \textit{loc. cit.} that

1) all of the above conjectures hold in the limit $k\to \infty$ (cf. \cite{BF1} for the exact formulation).

2) In \cite{BF1} a slightly weaker version of the above conjecture is proved in the case $k=1$;
the proof is based on the results of \cite{bfg}. A recent paper \cite{slofstra} proves it in full generality.

3) Again, a slightly weaker version of the above conjecture is proved in \cite{BF1} for $G=SL(N)$.
Let us mention the main ingredient of that proof.
Let $\grg$ be a simply laced simple finite-dimensional Lie algebra. Then by McKay correspondence
 one can associate with $\grg$ a finite subgroup $\Gam$ of $\SL(2,\CC)$.
Recall that H.~Nakajima (cf. e.g. \cite{nakajima})
gave a geometric construction of integrable $\grg_{\aff}$-modules of level $N$ using
certain moduli spaces which, roughly speaking, have to do
with vector bundles of rank $N$ on $\AA^2/\Gam$.
In particular, if $\grg=\fsl(k)$ it follows that

\smallskip
1) By H.~Nakajima the geometry of vector bundles of rank $n$ on $\AA^2/\Gam_k$ is related to integrable
modules over $\fsl(k)_{\aff}$ of level $N$.

2) By I.~Frenkel's suggestion  the geometry of vector bundles
of rank $N$ on $\AA^2/\Gam_k$ is related to integrable
modules over $\fsl(N)_{\aff}$ of level $k$.

\smallskip
\noindent
On the other hand, in the representation theory of affine Lie algebras
there is a well-known relation, due to I.~Frenkel,
between integrable
modules over $\fsl(k)_{\aff}$ of level $N$ and integrable
modules over $\fsl(N)_{\aff}$ of level $k$. This connection is called {\em level-rank duality}; one of its aspects
is discussed in \cite{Fr}. It turns out that combining the results of \cite{Fr} with the results of
\cite{nakajima} one can get a proof of a slightly weaker version of our main conjecture.

\medskip
\noindent
It is of course reasonable to ask why $G$-bundles on $\AA^2/\Gam_k$ have anything to do with the sought-for
affine Grassmannian of $G_{\aff}$. We don't have a satisfactory answer to this question, though some sort
of explanation (which would be too long to reproduce in the Introduction) is provided in \cite{BK}.
Also, E.~Witten produced an explanation of this phenomenon
in terms of 6-dimensional conformal field theory (cf. \cite{Wi}).

In~\cite{BF2} and \cite{BF3} other aspects of 1) are explored; in particular, an affine
analog of convolution of $G(\calO)$-equivariant
perverse sheaves on $\Gr_G$ and the analog
of the so called Beilinson-Drinfeld Grassmannian are discussed in \cite{BF2} and the analog of the Mirkovic-Vilonen fiber
functor is discussed in \cite{BF3}.
Most of the statements of \textit{loc. cit.} are still conjectural, but for the case of $G=SL(N)$ almost all of them
follow easily from the work of Nakajima \cite{Nak-Gr}.

\section{Eisenstein series for affine Kac-Moody groups}\label{aff-global}
In this Section we discuss unramified Eisenstein series for $G_{\aff}$ over a functional field;
we are going treat the subject in geometric way. Most of the results described below are adaptions
of the corresponding results of H.~Garland who considered the case of the global field $\QQ$); other
results are going to appear in \cite{BK-Eis}.

\subsection{The double quotient}As was explained in Section \ref{aut-finite} unramified automorphic
forms in the case of finite-dimensional $G$  are functions on the double quotient $G(\calO_{\AA_F})\backslash G(\AA_F)/G(F)$. We would like to introduce
an analogous quotient for the group $\hatG_{\aff}$ and then give some interesting examples of functions
on this quotient (given by an affine analog of Eisenstein series)
\footnote{Most of the discussion will go through for the group $G_{\aff}$ instead of $\hatG_{\aff}$, but for certain purposes which
go beyond the scope of this survey paper it seems more appropriate to work here with $\hatG_{\aff}$ here}.

Since $\hatG_{\aff}$ is a group ind-scheme over $\ZZ$, it makes sense to consider $\hatG_{\aff}(\AA_F)$.
However, it turns out that
this is not the right thing to consider. Instead, let us set
$$
\hatG_{\aff,\AA_F}= \{ (g_v\in \hatG_{\aff}(\calK_v))|\ g_v\in \hatG_{\aff}(\calO_v)\text{ for almost all $v$}\}.
$$
Then the double quotient on which we are going to produce some interesting functions is the quotient
$G_{\aff}(\calO_{\AA_F})\backslash G_{\aff,\AA_F}/G_{\aff}(F)$. To motivate the consideration of this double quotient
let us provide its geometric interpretation.
\subsection{The groupoid $\Bun_G(\hatS)$}
Let $S$ be a smooth
algebraic surface over a field $\kk$ and let $X\subset S$ be a smooth projective curve
in $S$ correspondiong to a sheaf of ideals $\calJ_X$. We let $S_n\subset S$ denote
the closed sub-scheme of $S$ corresponding to sheaf of ideals $\calJ_X^{n+1}$.
For each $n\geq m$ we have the embedding
$i_{mn}:S_m\to S_n$.

Let us denote by $\hatS$ the formal completion of $S$ along $X$ (it can be considered as either a formal scheme or an ind-scheme).
For an algebraic group $G$, we would like to consider the groupoid $\Bun_G(\hatS)$ of $G$-bundles on $\hatS$.
By the definition an object $\calF$ of $\Bun_G(\hatS)$ consists of
$G$-bundles $\calF_n$ on each $S_n$ together with a compatible system of isomorphisms
$\calF_n|_{S_m}\simeq \calF_m$ (and the notion of isomorphism of such objects is clear).
\subsection{$\Bun_G(\hatS^0)$}
We start our definition of $G$-bundles on $\hatS^0=\hatS\backslash X$ with the case $G=GL(n)$. In other words, we are going to
define the category $\Vect(\hatS^0)$ of vector bundles on $\hatS^0$.

Let us denote by $\Vect(\hatS)$ the category of locally free coherent sheaves on $\hatS$;
$\bullet$ Objects of $\Vect(\hatS^0)$ are objects of $\Vect(\hatS)$.

$\bullet$ Given two objects $\calE_1,\calE_2$ of $\Vect(\hatS)$ we define
$$
\Hom_{\Vect(\hatS^0)}(\calE_1,\calE_2)=\bigcup\limits_{n\geq 0}\Hom_{\Vect(\hatS)}
(\calE_1,\calE_2(nX)).
$$

Let us now assume that $G$ is simply connected.
Since $\Vect(\hatS^0)$ is a tensor category we can define a groupoid $\Bun_G'(\hatS^0)$
of exact tensor functors from  $\Rep(G)$ to $\Vect(\hatS^0)$.  It is clear that we have
a functor $r: \Bun_G(S)\to \Bun_G'(\hatS^0)$ and we  define $\Bun_G(\hatS^0)$ to be the full subcategory
of $\Bun_G'(\hatS^0)$ consisting of objects which lie in the image of $r$.

We conjecture that $\Bun_G(\hatS^0)=\Bun_G'(\hatS^0)$.

Since we assume that $G$ is simply connected by  Theorem 11.5 of \cite{Starr} this conjecture is true if we replace $\hatS$ by $S$.

\subsection{Bundles with respect to $G_{\aff}'$}
We would like to define the groupoid $\hatG_{\aff}$-bundles on $X$; let us first do it for the group $\hatG'_{\aff}$.
 Since we are given a homomorphism
$\hatG'_{\aff}\to \GG_m$, if such a notion makes sense, then we should have a functor
$$
\iota:\Bun_{G'_{\aff}}\to \Pic(X).
$$
So, in order to describe $\Bun_{G'_{\aff}}$, it is enough to describe the groupoid $\eta^{-1}(\calL)$ for each
$\calL\in \Pic(X)$ in a way that the assignment $\calL\mapsto \eta^{-1}(\calL)$ is functorial with
respect to isomorphisms in $\Pic(X)$. We set
$$
\iota^{-1}(\calL)=\Bun_G(\hatS_{\calL}^0)
$$
where $S$ is the total space of $\calL$.

\begin{proposition}\label{bun=adeles-affine}
We have an equivalence of groupoids
$$
\Bun_{\hatG'_{\aff}}\simeq \hatG'_{\aff}(\calO(\AA_F))\backslash \hatG'_{\aff,\AA_F}/\hatG'_{\aff}(F).
$$
\end{proposition}
Let us now turn on the central extension. Let $X\subset S$ be as in the beginning of this Section (in particular, we do not assume
anything about the self-intersection of $X$).

The central extension
$$
1\to \GG_m\to \hatG_{\aff}\to \hatG_{\aff}'\to 1
$$
gives rise to a $\Pic(X)$-torsor $\tBun_G(\hatS^0)$ over $\Bun_G(\hatS^0)$ for any $X\subset S$ as above.
Denote as as before by $\Bun_{\hatG_{\aff}}$ the union of $\tBun_G(\hatS_{\calL}^0)$ where $\calL$ runs over all
elements of $\Pic(X)$. Then the analog of
Proposition \ref{bun=adeles-affine} reads as follows:
\begin{equation}\label{bun=adeles-affinee}
\Bun_{\hatG_{\aff}}\simeq \hatG_{\aff}(\calO(\AA_F))\backslash \hatG_{\aff,\AA_F}/\hatG_{\aff}(F).
\end{equation}
\subsection{Eisenstein series}We now fix $\calL\in \Pic(X)$ and assume that $\deg(\calL)<0$.
For a parabolic subgroup $P\subset G$, let $\Bun_{G,P}(\hatS_{\calL})$ be the groupoid of
$G$-bundles on $\hatS_{\calL}$ endowed with a $P$-structure on $X\subset \hatS_{\calL}$. Note that
$\Bun_{G,G}(\hatS_{\calL}$) is just $\Bun_G(\hatS_{\calL})$ (in the future we shall use the following
convention: for any symbol of the form $?_{G,P}$, defined for all parabolic subgroups of $G$, we shall
write $?_G$ instead of $?_{G,G}$). Consider
diagram
\begin{equation}\label{aff-diag}
\begin{CD}
\Bun_{G,P}(\hatS_{\calL})\x\Pic(X)@>\eta_{G,P,\calL}>>\Bun_M(X)\x\Pic(X)\\
@V\pi_{G,P,\calL} VV\\
\tBun_G(\hatS^0_{\calL})
\end{CD}
\end{equation}
We shall mostly be interested in the cases $P=G$ and $P=B$ (a Borel subgroup of $G$);
for $P=G$ we shall just write $\Eis_{G,\calL}$ instead of $\Eis_{G,G,\calL}$.

Similarly to Subsection \ref{e-finite} we let $\rho_P^{\aff}:\ZZ\x\Lam_M\to \ZZ$ be the homomorphism which such that
$\rho_{P,\aff}(\alp)=0$ for every coroot $\alp$ of $M$ and $\rho_{P,\aff}(\alp_i)=1$ for every simple coroot
of $\hatG_{\aff}$ not lying in $M$. We also denote by $\rho_{P,\aff}$ the corresponding function
$\Pic(X)\x\Bun_M(X)\to \ZZ$. For a complex-valued function $f$ on $\Bun_M(X)\x\Pic(X)$ we would like to set
\begin{equation}\label{aff-eis}
\Eis_{G,P,\calL}(f)=(\pi_{G,P,\calL})_!\eta_{G,P,\calL}^*(q^{\rho_{P,\aff}}f).
\end{equation}
A priori, it is not clear what sense it makes, since the RHS of (\ref{aff-eis}) might be an infinite sum.
We say that
$\Eis_{G,P,\calL}(f)$ is well-defined if  it value at every point is
given by a finite sum
Then we have
\begin{theorem}\label{aff-eis-one}
\begin{enumerate}
\item Assume that $f$ has finite support. Then $\Eis_{G,P,\calL}(f)$ is well-defined (i.e. its value at every point is
given by a finite sum). If $P=G$ then the same is true for any $f$ such that the image of $\supp(f)$ in $\Pic(X)$ is finite.
\item
Assume that $P=B$; then $\Bun_M(X)\x\Pic(X)=\Bun_T(X)\x\Pic(X)$. Let $\deg:\Bun_T(X)\x\Pic(X)\to\Lam\oplus\ZZ$ denote
the natural degree map and let $f_{\grs}=\chi\cdot q^{\la\deg,\grs\ra}$, where $\chi$ is a unitary character of $\Bun_T\x\Pic(X)$
and $\grs\in\grt^*\oplus\CC$. Then the series
$\Eis_{G,B,\calL}(f_{\grs})$ is absolutely convergent when $\Re \grs > \rho_{\aff}^{\vee}$.
Moreover, it has a meromorphic continuation to the domain $\Re \grs>0$.
\item
Assume that $P=G$. Then $\Eis_{G,P,\calL}(f)$ is well-defined if $f$ has finite support modulo $\Pic(X)$. In particular,
this is true for any cuspidal $f$.
\end{enumerate}
\end{theorem}

The first statement is due to Kapranov \cite{Kap-E}. The second statement is due to Garland \cite{Gar1} and the third
statement will appear in \cite{BK-Eis}.

Given $P$ as above one can consider the parahoric subgroup $\calP_P\subset \hatG_{\aff}$ defined as
$\GG_m\x G[[t]]_P\rtimes \GG_m$, where $G[[t]]_P$ is the preimage of $P$ under the natural map $G[[t]]\to G$.
It should be thought of as a parabolic subgroup of the affine Kac-Moody group $\hatG_{\aff}$ and the corresponding
Levi factor is $\GG_m\x M\x\GG_m$ where $M$ is the Levi factor of $P$.
The operators $\Eis_{G,P,\calL}$ should be thought of as Eisenstein series $\Eis_{\hatG_{\aff},\calP_P}$ (restricted
to $\iota^{-1}(\calL)$). If $B\subset G$ is a Borel subgroup, then $\calP_B$ should be thought of as a Borel
subgroup of $\hatG_{\aff}$ (later, we are going to consider a different type of Borel subgroups of $\hatG_{\aff}$.
In fact, it is well-known that the group $\hatG_{\aff}$ has more general parabolic subgroups containing
$\calP_B$ than those  of the form $\calP_P$ for some $P\subset G$; we shall call such subgroups \textit{parabolic subgroups
of $\hatG_{\aff}$
of positive type} (since later we are going to consider other parabolic subgroups). The above Eisenstein series
can be defined for any parabolic subgroup $\calP\subset \hatG_{\aff}$ of positive type and the analog of Theorem \ref{aff-eis-one}
holds for any such $\calP$.

\subsection{Constant term}\label{constant-one}
We now want to define the operator of constant term, acting from functions on $\tBun_G(\hatS^0_{\calL})$ to functions on
$\Pic(X)\x \Bun_M(X)$. In principle, we would like to set
$c_{G,P,\calL}(f)=q^{\rho_{P,\aff}}(\eta_{G,P,\calL})_!\pi_{G,P,\calL}^*(f)$. However, formally there is a problem with this
definition. Namely, the morphism $\eta_{G,P,\calL}$ is not representable and the corresponding automorphism groups
(whose sizes enter the definition of $(\eta_{G,P,\calL})_!$ as in (\ref{pshrieck})) are infinite. So, we need to apply
certain "renormalization" procedure in order to define $(\eta_{G,P,\calL})_!$. Let us explain this procedure for
$P=G$(the general case is simila). Let $\calF\in \Bun_G(\hatS_{\calL})$.
The, for any $\calM\in \Pic(X)$, the automorphism group of $(\calM,\calF)$ as a point in
$\eta_{G,\calL}^{-1}(\calM,\calF|_X)$ is $\Aut^0(\calF)$ which consists of those automorphisms of $\calF$ which are
equal to identity when restricted to $X$. This is an infinite group which is equal to the set of $\FF_q$-points of
some pro-unipotent algebraic group and thus formally its "number of points" should be the same as the number of
points in the corresponding Lie algebra, which is equal to $H^0(X,p_*(\grg_{\calF})(-X))$, where $p:\hatS_{\calL}\to X$ denotes the
natural morphism and $\grg_{\calF}$ is the vector bundle associated to $\calF$ by means of the adjoint representation of $G$.
The space $H^0(X,p_*(\grg_{\calF})(-X))$ is infinite-dimensional, but the corresponding space $H^1(X,p_*(\grg_{\calF})(-X))$
is actually finite-dimensional (here we use the fact that $\deg(\calL)<0$) and we formally set
$$
\#\Aut^0(\calF)=q^{\dim H^1(X,p_*(\grg_{\calF})(-X))}.
$$
Informally, this definition is reasonable because by the Riemann-Roch theorem if the Euler characteristic
$\dim H^0(X,p_*(\grg_{\calF})(-X))-\dim H^1(X,p_*(\grg_{\calF})(-X))$ made sense, it would have been independent of
$\calF$. Given this convention, we can now define $(\eta_{G,P,\calL})_!$ and for a function $f$ on
$\tBun_G(\hatS^0_{\calL})$ we set
$$
c_{G,P,\calL}(f)=q^{-\rho_{P,\aff}}(\eta_{G,P,\calL})_!\pi_{G,P,\calL}^*(f).
$$
In fact, a slight variation gives
a definition of $c_{\hatG_{\aff},\calP}$ for any parabolic subgroup $\calP$ of $\hatG_{\aff}$ of positive type.
\footnote{Note that the factor $q^{\rho_{P,\aff}}$ is changed to $q^{-\rho_{P,\aff}}$; this is due to the fact
a factor of $q^{2\rho_{P,\aff}}$ is hidden
in the renormalization procedure described above}.
For $P=B$ it is easy to see that this definition coincides with that of \cite{Gar1}. It is easy to see that
$c_{G,P,\calL}(f)$ makes sense for \textit{any} function $f$ on $\tBun_G(\hatS^0_{\calL})$, since any fiber of the map
$\Bun_G(\hatS_{\calL})\to \Bun_G(X)$ has finitely many isomorphism classes of objects.
\subsection{More on Eisenstein series and constant term}Given a parabolic subgroup $P$ of $G$ one consider the subgroup
$\calQ_P$ of $\hatG_{\aff}$ defined
$\GG_m\x G[t^{-1}]_P\rtimes \GG_m$ where $G[t^{-1}]_P$ is the preimage of $P$ under the map
$G[t^{-1}]\to G$ obtained by setting $t^{-1}=0$. We shall call such subgroups \textit{parabolic subgroups
of negative type} (as before there exist more general parabolic subgroups of negative type, but we shall
not consider them in this paper). One can talk about constant term and Eisenstein series for
parabolic subgroups of negative type. Geometrically, these are defined by means of the following analog of
(\ref{aff-diag}).
\begin{equation}\label{aff-diag-}
\begin{CD}
\Pic(X)\x\Bun_{G,P}(S_{\calL^{-1}})@>\eta^-_{G,P,\calL}>>\Pic(X)\x\Bun_M(X)\\
@V\pi^-_{G,P,\calL} VV\\
\tBun_G(\hatS^0_{\calL})
\end{CD}
\end{equation}
Here $S_{\calL^{-1}}$ denotes the total space of $\calL^{-1}$ and $\Bun_{G,P}(S_{\calL^{-1}})$
denotes the groupoid of $G$-bundles on $S_{\calL^{-1}}$ endowed with a $P$-structure on $X\subset S_{\calL^{-1}}$.
The existence of the map $\Bun_{G,P}(S_{\calL^{-1}})\to \Bun_G(\hatS^0_{\calL})$ is clear; its lift to a map
$\Pic(X)\x\Bun_{G,P}(S_{\calL^{-1}})\to\tBun_G(\hatS^0_{\calL})$ is not difficult to define, but we shall not do it in this paper.
Thus one can define operator $c_{G,P,\calL}^-$.
In this case no "renormalization" as in Subsection \ref{constant-one} is needed, but the fibers of the map $\eta^-_{G,P,\calL}$
are infinite and thus a priori $c_{G,P,\calL}^-$ can only be applied to functions with finite support -- otherwise
one needs to check convergence (this is very different from the finite-dimensional case).

One can also define the corresponding Eisenstein series operator $\Eis_{G,P,\calL}^-$; this is an operator from
functions with finite support on $\Pic(X)\x\Bun_G(X)$ to functions on $\tBun_G(\hatS^0_{\calL})$. However,
its definition is not straightforward,
since the composition $(\pi^-_{G,P,\calL})_!(\eta^-_{G,P,\calL})^*(f)$ is given by an infinite sum even when $f$
has finite support and we must use a renormalization procedure similar to that from Subsection \ref{constant-one}.
The formal definition will appear in \cite{BK-Eis}.

We now want to present the analog Theorem \ref{gk}.
Let $\grn_{P,\aff}$ be the nilpotent radical of Lie$\calP_P$ and let $\grn_{P,\aff,-}$ be the nilpotent
radical of Lie$\calQ_P$. We can also consider the corresponding Langlands dual subalgebras
in $\grg_{\aff}^{\vee}$.
\begin{theorem}\label{gk-aff}
Assume that $f\in \CC(\Pic(X)\x\Bun_M(X))$ is a cuspidal Hecke eigen-function. Then
\begin{enumerate}
\item
$$
\begin{aligned}
&c_{G,P,\calL}\circ \Eis_{G,P,\calL}(f)=\\
\sum\limits_{w\in W_{\aff}/W(M)} &
q^{-(g-1)\dim \grn_{P,\aff}^{\vee}\cap (\grn_{P_-,\aff,-}^{\vee})^w}f^w
\frac{L(f, \grn_{P,\aff}^{\vee}\cap (\grn_{P_-,\aff,-}^{\vee})^w, 0)}
{L(f, \grn_{P,\aff}^{\vee}\cap (\grn_{P_-,\aff,-}^{\vee})^w, 1)}.
\end{aligned}
$$
Note that the factor in front of $f^w$ is slightly different from that in Theorem \ref{gk}(1); this has to do with the
renormalization procedure of Subsection \ref{constant-one}.
\item
Assume that $G$ is simply laced and let $f=(\chi,\phi)$ where $\chi$ is a character of $\Pic(X)$ and $\phi$
is a cuspidal Hecke eigen-function on $\Bun_M(X)$. Then
$$
\begin{aligned}
&c^-_{G,P,\calL}\circ \Eis_{G,P,\calL}(f)=
\prod\limits_{i=1}^r\prod\limits_{j=1}^{\infty} \frac{L(\chi^j,d_i)}{L(\chi^j, d_i-1)}\x\\
\sum\limits_{w\in W_{\aff}/W(M)} &
q^{(g-1)\dim \grn_{P,\aff}^{\vee}\cap (\grn_{P_-,\aff,-}^{\vee})^w}f^w
\frac{L(f, \grn_{P,\aff}^{\vee}\cap (\grn_{P,\aff}^{\vee})^w, 0)}
{L(f, \grn_{P,\aff}^{\vee}\cap (\grn_{P,\aff}^{\vee})^w, 1)}.
\end{aligned}
$$
A similar statement holds for the composition $c_{G,P,\calL}\circ \Eis_{G,P,\calL}^-$.
\end{enumerate}
\end{theorem}

A few remarks are in order about the formulation of Theorem \ref{gk-aff}. The first statement is proved in a similar manner
to the finite-dimensional case; in the case $P=B$ it appears in \cite{Gar1} (we also refer the reader to \cite{Gar1} for the
discussion of convergence of the right hand side of Theorem \ref{gk-aff}(1)).
In the second statement the product of abelian $L$-functions in front of the
sum comes from the formula (\ref{delta}) for the "correction term" in the affine Gindikin-Karpelevich formula;
a variant of this statement holds for non-simply laced groups as well (using the formula for the correction term described
in \cite{BFK}). Let us now look at the case $P=G$. Then in the right hand side of Theorem \ref{gk-aff}(2) only one
term (corresponding to $w=1$) is left and it is equal to $f$ multiplied by certain infinite product of ratios of $L$-functions.
Moreover, let us consider the case when $\chi$ above is equal to $q^{s\deg}$. Then the product
$$
\prod\limits_{i=1}^r\prod\limits_{j=1}^{\infty} \frac{L(\chi^j,d_i)}{L(\chi^j, d_i-1)}=
\prod\limits_{i=1}^r\prod\limits_{j=1}^{\infty} \frac{\zeta(js+d_i)}{\zeta(js+d_i-1)},
$$
which is a meromorphic function of $s$ when $\Re(s)>0$. Also,
$$
\frac{L(f, \grn_{P,\aff}^{\vee}, 0)}
{L(f, \grn_{P,\aff}^{\vee}, 1)}=
\prod\limits_{j=1}^{\infty} \frac{L(\phi,\rho,js)}{L(\phi,\rho,js+1)},
$$
where $\rho$ is the adjoint representation of $G^{\vee}$. The above product is absolutely convergent for $\Re (s) >1$ and the standard
conjectures about automorphic $L$-functions imply that it should have a meromorphic continuation to the domain $\Re(s)>0$. We do not
know how to prove this at the moment, but we expect this observation to be a useful tool in proving that $L(\phi,\rho,s)$
has a meromorphic continuation.

It is also important to note that we are not claiming anything about the composition $c_{G,P,\calL}^-\circ\Eis_{G,P,\calL}^-$
since we don't know how to make sense of it.
\subsection{Affine Tamagawa number formula}
To conclude this Section, we would like to present an analog of Theorem \ref{tam-finite} in the affine case. First, one needs
to define a measure on $\tBun_G(\hatS^0_{\calL})$ analogous to the one defined by (\ref{measure}). Naively, one can try to define
in it by the same formula as (\ref{measure}) on $\Bun_G(\hatS^0_{\calL})$ rather than on $\tBun_G(\hatS^0_{\calL})$;
however, one quickly discovers that the automorphism groups which appear in (\ref{measure}) are infinite, and to define the measure
one needs to perform again a renormalization procedure similar to the one in Subsection \ref{constant-one}. We are not going to
give the details of that procedure here, but let is just mention that it is a little more involved than in Subsection \ref{constant-one}
and, in particular, the resulting measure makes sense only on $\tBun_G(\hatS^0_{\calL})$; moreover, with respect to the natural
$\Pic(X)$-action on it changes according to the character
$\calM\mapsto q^{-2h^{\vee}\deg(\calM)}$. This phenomenon has been discovered by H.~Garland in \cite{Gar1}, who gave a group-theoretic
definition of this measure.

Since the above measure changes according to a non-trivial character of $\Pic(X)$, we can't integrate the function 1 with
respect to this measure. Instead, we are going to define the volume of $\tBun_G(\hatS^0_{\calL})$ by recalling the interpretation
of the volume of $\Bun_G(X)$ given by Theorem \ref{finite-residue}. Namely, for $\grs\in\CC\oplus\grt^*$
we set $E(\grs)=\Eis_{G,B,\calL}(f_{\grs})$
where $f_{\grs}$ is equal to $\la (n,\gam),\grs\ra$ on $\Pic^n(X)\x\Bun_T^{\gam}(X)$. Thus we can define
$\grr_{\hatG_{\aff}}$ as the residue of $c_{G,B,\calL}^-(E(\grs))$ at $\grs=\rho_{\aff}$ multiplied by
by the function $q^{\la ?,\rho_{\aff}\ra}$ (this is a constant function on $\Pic(X)\x\Bun_T(X)$
which priori, $\grr_{\hatG_{\aff}}$ might depend on $\calL$ but it is easy to deduce from
Theorem \ref{gk-aff} that it does not). Then we define the volume of $\tBun_G(\hatS^0_{\calL})$ as
$$
\frac{(\Res_{s=1}\zeta(s))^{r+1}}{\grr_{\hatG_{\aff}}}.
$$

With this definition we now have the following
\begin{theorem}\label{aff-tamagawa}
Assume that $G$ is simply laced and simple. Let us assume that the exponents $d_1,\cdots, d_r$ are numbered
in such a way that $d_1=2$. Then the volume of $\tBun_G(\hatS^0_{\calL})$ is equal to
\begin{equation}\label{aff-tam-form}
\frac{\prod\limits_{i=1}^r \zeta(d_i)}{\prod\limits_{i=2}^r \zeta(d_i-1)}.
\end{equation}
\end{theorem}
Let us make
two remarks about (\ref{aff-tam-form}). First, the absence of the factor corresponding to $d_1$ in the denominator has
to do with the fact that we are working not with the group $\hatG'_{\aff}$ but with its central extension $\hatG_{\aff}$.
Second, usually most terms in (\ref{aff-tam-form}) cancel out (e.g. for $G=SL(n)$ only $\zeta(n)$ survives). However, writing
the answer as in (\ref{aff-tam-form}) is still instructive, since
one can give an explanation why the above answer is natural in the spirit of Theorem \ref{at-bott} (however, we don't
know how to formulate precisely a "cohomological" statement that would imply Theorem \ref{aff-tamagawa}).

\section{Further questions and constructions}\label{further}
In this Section we mention some related works and formulate several possible directions of future research.
\subsection{Existence of cuspidal representations and automorphic forms}
It is interesting whether there exists a reasonable notion of cuspidal automorphic forms for $\hatG_{\aff}$.
We can try to say that a function $f$ on $\tBun_G(\hatS^0_{\calL})$ is cuspidal if $c_{\hatG_{\aff},\calP}(f)=0$
for every parahoric subgroup $\calP$ of $\hatG_{\aff}$. Here the words "for every parahoric" can in principle be interpreted
in several different ways. Namely, we can require that this holds for all parahoric subgroups of positive type, or for
all parahoric of negative type, or both. We conjecture that these conditions are in fact equivalent. More precisely,
we conjecture that if $\calP$ and $\calQ$ are opposite parahoric subgroups then
$c_{\hatG_{\aff},\calP}(f)=0$ implies that $c_{\hatG_{\aff},\calQ}=0$.

If the above conjecture is true, then we get an unambiguous definition of cuspidal functions. However, it is not at all
clear whether non-zero cuspidal functions exist. It would be very interesting to understand if they do exist and how to describe
the space of cuspidal functions.

The above questions make sense locally.
As was mentioned in the introduction, in the paper \cite{GK1}-\cite{GK2} the theory of representations
of the group $\hatG_{\aff}(\calK)$ (for $\calK$ being a local non-archmedian field) was developed.
One can define the notion of a cuspidal representation in this framework and it would be very interesting to understand whether
$\hatG_{\aff}$ has any irreducible cuspidal representations.

\subsection{Weil representation and theta-correspondence}
In \cite{Zhu} Y.~Zhu generalized the notion of Weil representation of the double cover $\widetilde{\Sp}(2n,\calK)$
of the symplectic group $\Sp(2n,\calK)$ over a local field $\calK$ (archimedian or not) to its affine analog;
moreover, the affine analog of the corresponding automorphic representation is also constructed in \textit{loc. cit.}

One of the main applications of the usual Weil representation to the representation theory of $p$-adic groups and to automorphic
forms is the construction of the so called \textit{theta-correspondence} (cf. \cite{pras} for a survey on theta-correspondence).
One of the most interesting features of the theta-correspondence is that it provides a fairly explicit tool for producing
examples of cuspidal automorphic forms.

It would be very interesting to develop an affine analog of theta-correspondence; in
particular, it might give rise to a construction of cuspidal automorphic forms for
the group $\hatG_{\aff}$ (for some particular choices of $G$). The full theory of theta-correspondence has not yet been developed in the affine case.
However, it was shown in \cite{Gar-Zh1} and \cite{Gar-Zh2} that one of the main tools used in the theory
of (global) theta-correspondence -- the so called \textit{Siegel-Weil formula} -- does have an analogue in the affine case.

\subsection{Meromorphic continuation and functional equation of Eisenstein series}As was claimed above, the Eisenstein series
$\Eis_{G,\calL}(f)$ is
convergent for any cuspidal function $f$ on $\Pic(X)\x\Bun_G(X)$. This fact is not true for
$\Eis^-_{G,\calL}$. Let $f$ be of the form $\chi\x\phi$ where $\phi$ is a cuspidal function on
$\Bun_G(X)$ and $\chi$ is a unitary character of $\Pic(X)$. For every $s\in \CC$ set
$$
f_s(\calM,\calF)=\chi(\calM)\phi(\calF)q^{s\deg(\calM)}.
$$
Then one can show that $\Eis^-_{G,\calL}(f_s)$ is absolutely convergent for $\Re(s)>h^{\vee}$ where $h^{\vee}$ is the dual
Coxeter number of $G$. We conjecture that in fact $\Eis^-_{G,\calL}(f_s)$ has a meromorphic continuation to
the domain $\Re(s)>0$.

Moreover, we conjecture that $\Eis^-_{G,\calL}(f_s)$ is proportional to $\Eis_{G,\calL}(f_s)$ (the coefficient of proportionality
can then be expressed as a ratio of an infinite product of $L$-functions as follows from

\subsection{Towards Kazhdan-Lusztig theory for DAHA}
Let us recall the algebra $\bfH^+_{\aff}$ defined in Subsection \ref{aff-iwahori}. Moreover, as was remarked after
Theorem \ref{aff-iw-theorem}, when we choose a local non-archimedian field $\calK$ the specialization of  $q$ to
the number of elements in its residue field acquires a natural basis. We conjecture that this basis comes from a
$\CC[q,q^{-1}]$-basis of $\bfH^+_{\aff}$; it should be thought of as an analog of the "standard basis" of a finite
of affine Hecke algebra. It would be interesting to define an analogue of Kazhdan-Lusztig basis of $\bfH^+_{\aff}$.
We don't know how to attack this problem algebraically; however, it should be possible to attack it geometrically using
appropriate generalization of the constructions discussed in Subsection \ref{aff-geom-satake}.

\end{document}